        \documentstyle[11pt,amsmath,amsfonts,euscript]{article}
        \newcommand{\al}{\alpha}
        
        \newcommand{\del}{\delta}
        \newcommand{\eps}{\epsilon}
        \newcommand{\lam}{\lambda}
        \newcommand{\sig}{\sigma}
        
        \newcommand{\vth}{\vartheta}
        
        \newcommand{\Del}{{\mathit{\Delta}}}
        \newcommand{\Gam}{{\mathit{\Gamma}}}
        \newcommand{\Lam}{{\mathit{\Lambda}}}

        \newcommand{\LL}{{\cal L}}
        \newcommand{\OO}{{\cal O}}
        
        \newcommand{\GG}{{\EuScript{G}}}
        \newcommand{\MM}{{\EuScript{M}}}

        \newcommand{\g}{{\mathfrak g}}
        \newcommand{\gt}{{\mathfrak t}}
        
        \newcommand{\ZA}{{\mathbb{Z}}}
        \newcommand{\RE}{{\mathbb{R}}}
        \newcommand{\CO}{{\mathbb{C}}}

        \newcounter{sect}
        \newcounter{subsect}
        \newcommand{\sect}[1]{\vspace{2ex}
                \addtocounter{sect}{1}\setcounter{subsect}{0}
                \begin{flushleft}
                {{\large\bf \arabic{sect}. {#1}}}
                \end{flushleft}
                \setcounter{thm}{0}
                \setcounter{equation}{0}
                \def\theequation{\arabic{sect}.\arabic{equation}}
                \def\thefigure{\arabic{sect}.\arabic{figure}}}
        
        \newtheorem{thm}{Theorem}[sect]
        \newtheorem{prop}[thm]{Proposition}
        \newtheorem{lemma}[thm]{Lemma}

        \newcommand{\be}{\begin{equation}}
        \newcommand{\ee}{\end{equation}}
        \newcommand{\bea}{\begin{eqnarray}}
        \newcommand{\eea}{\end{eqnarray}}
        \newcommand{\nno}{\nonumber \\}

        \newcommand{\ii}{\sqrt{-1}\,}

        \newcommand{\e}[1]{e^{{#1}}}

	\newcommand{\Sum}[1]{\underset{{}^{#1}}{\sum}}
        \newcommand{\bra}{\langle}
        \newcommand{\ket}{\rangle}
        
        \newcommand{\vol}{{\rm vol}\/}
        \newcommand{\ad}{{\rm ad}\/}
        \newcommand{\id}{{\rm id}\/}
        \newcommand{\Hom}{{\rm Hom}\/}
        \newcommand{\spanz}{{\rm span}_\ZA}
        \newcommand{\spanc}{{\rm span}_\CO}
        \newcommand{\loong}{_{\rm long}}
        \newcommand{\short}{_{\rm short}}
        
        \newcommand{\set}[2]{\{{#1}|{#2}\}}
	\newcommand{\sqrtn}{\sqrt{n}}
        \newcommand{\LA}{{}^L\!A}
        \newcommand{\LLL}{{}^L\!\LL}
        \newcommand{\LLam}{{}^L\!\Lam}

\begin{document}
\begin{flushright}
{\tt arXiv:yymm.nnnn [math.RT]}
\end{flushright}

\vspace{5pt}

        \begin{center}
{\Large\bf Miniscule representations, Gauss sum and modular invariance}
\footnote{Talk at the 4th International Congress of Chinese Mathematicians,
Hangzhou, December 2007}\\
        \vspace{4ex}
        {\large\rm Siye Wu}

	\vspace{2ex}

{\small {\em Department of Mathematics, University of Hong Kong, Pokfulam,
 Hong Kong, China}\footnote{Current address. E-mail: {\tt swu@maths.hku.hk}}

and

       {\em Department of Mathematics, University of Colorado, Boulder,
          CO 80309-0395, USA}}
        \end{center}

        \vspace{2ex}

        \begin{quote}
{\small 
\noindent {\bf Abstract.}
After explaining the concepts of Langlands dual and miniscule 
representations, we define an analog of the Gauss sum for any compact,
simple Lie group with a simply laced Lie algebra.
We then show a reciprocity property when a Lie group is exchanged with
its Langlands dual.
We also explore the relation with theta functions and modular transformations.
In the non-simply laced case, we construct a unitary representation of the
Hecke group which involves interesting new phase factors.


}
	\end{quote}

\vskip 1cm

\sect{Introduction}

Duality has played an important role in physics and mathematics.
One of the earliest examples is the electric-magnetic duality or $S$-duality 
[\ref{GNO}, \ref{MO}], which is exact in the $N=4$ supersymmetric gauge theory.
This leads to the celebrated conjecture of Vafa and Witten [\ref{VW}] on
the Euler number of the moduli space of anti-self-dual connections.
In this paper, we study some finite dimensional unitary representations of
the modular group motivated from the $S$-duality transformations in 
[\ref{VW}] and their consequences.
We also construct similar representations of the Hecke groups.
We leave the interplay with four manifolds to subsequent work.

The paper is organized as follows.
In Section~2, we review the concepts of Langlands dual and miniscule
representations.
In Section~3, we study the simply laced case.
We define an analog of the Gauss sum for a compact simple Lie group 
using miniscule representations.
We then show a reciprocity property when the Lie group is exchanged with
its Langlands dual.
We also explore the relation with theta functions and modular transformations.
Section~4 is on the non-simply laced case.
We start with the theta functions and derive an analog of the
modular transformations for the Hecke group.
This enables us to construct a finite dimensional unitary representation
of the Hecke group which involves new phase factors that will have
interesting implications on $S$-duality.

\sect{Langlands dual and miniscule representations}

To explain the concept of dual groups, we start with the Abelian case.
Suppose $T$ is a compact real torus with Lie algebra $\gt$.
Then $T=\gt/2\pi\ii\ell$, where $\ell\subset\ii\gt$ is a lattice
of full rank. 
We note that $\ell\cong\Hom(U(1),T)$.
On the other hand, the (one-dimensional) irreducible representations of $T$
are classified by the dual lattice $\ell^*\subset\ii\gt^*$.
The dual torus is $T^*=\gt^*/2\pi\ii\ell^*$, for which the roles of $\ell$
and $\ell^*$ are reversed.
That is, $\ell^*\cong\Hom(U(1),T^*)$ whereas $\ell$ classifies the irreducible
representations of $T^*$.

Now let $G$ be a simple, compact and connected Lie group.
The Langlands dual ${}^LG$ of $G$ is characterised by the property that
the irreducible representations of ${}^LG$ are in one-to-one correspondence
with the homomorphisms from $U(1)$ to $G$ modulo the conjugation in $G$.
We will describe a construction of ${}^LG$ below.
Some examples of $G$ and its dual ${}^LG$ are listed as follows.
\begin{center}
\begin{tabular}{|c|cccccc|}  \hline
$G$ & & $SU(n)$ & $Spin(2n)$ & $Sp(n)$ & $Spin(2n+1)$ & $E_8$ \\ \hline
${}^LG$ & & $SU(n)/\ZA_n$ & $SO(2n)/\ZA_2$ & $SO(2n+1)$ & $Sp(n)/\ZA_2$ 
& $E_8$ \\ \hline
\end{tabular}
\end{center}

Let $T$ be a maximal torus of $G$ and $\g$, $\gt$, the Lie algebras of
$G$, $T$, respectively.
Recall that the root system $\Del\subset\ii\gt^*$ is the set of roots and
that the root lattice $\Lam=\spanz\Del\subset\ii\gt^*$ is the lattice
generated by $\Del$.
For any $\al\in\Del$, the corresponding coroot $\check\al\in\gt$ is the
vector such that $\bra\beta,\check\al\ket=2(\beta,\al)/(\al,\al)$ for
any $\beta\in\Del$.
Here $\bra\cdot,\cdot\ket$ is the pairing between $\ii\gt$ and it dual
space $\ii\gt^*$ while $(\cdot,\cdot)$ is an inner product on $\ii\gt^*$
invariant under the Weyl group; such an inner product is unique up to a 
multiple of positive scalars (since $\g$ is simple) and the coroot 
$\check\al$ does not depend on its choice.
The set $\check\Del=\set{\check\al}{\al\in\Del}$ is also a root system.
Let $\check\Lam=\spanz\check\Lam\subset\ii\gt$ be the coroot lattice.
Then the weight lattice is $\check\Lam^*\subset\ii\gt^*$ and the coweight
lattice is $\Lam^*\subset\ii\gt$.
Let $\ell\subset\ii\gt$ be the lattice such that $T=\gt/2\pi\ii\ell$.
Then we have the inclusions [\ref{B9}, \S 4.9]
\be
\check\Lam\subset\ell\subset\Lam^*\subset\ii\gt, \quad
\Lam\subset\ell^*\subset\check\Lam^*\subset\ii\gt^*.
\ee

The Lie algebra $\g$ uniquely determines the universal covering group
$\tilde G$ of $G$, which is simply connected and whose centre is
$Z(\tilde G)\cong\Lam^*/\check\Lam$.
Here is a table of $Z(\tilde G)$.
\begin{center}
\begin{tabular}{|c|ccccccc|}\hline
$\g$ & $A_r$ & $B_r$ & $C_r$ & $D_r$ & $E_{r=6,7,8}$ & $F_4$ & $G_2$ \\ \hline
$Z(\tilde G)$ & $\ZA_{r+1}$ & $\ZA_2$ & $\ZA_2$ & 
$\!\underset{(r \text{ odd})}{\ZA_4}
\underset{(r \text{ even})}{\ZA_2\oplus\ZA_2}\!$ & $\ZA_{9-r}$ & $1$ & $1$ \\
\hline
\end{tabular}
\end{center}
At the opposite extreme of $\tilde G$, the adjoint group 
$G_\ad=\tilde G/Z(\tilde G)$ has the same Lie algebra $\g$, 
but has $Z(G_\ad)=1$ and $\pi_1(G_\ad)=Z(\tilde G)$.
Notice that $\ell=\check\Lam$ for $G=\tilde G$ and $\ell=\Lam^*$ for $G=G_\ad$.
For a general compact group $G$ with Lie algebra $\g$, we have
$Z(G)=\Lam^*/\ell$ and $\pi_1(G)=\ell/\check\Lam$.

The construction of the Langlands dual ${}^LG$ is achieved in two steps.
First, the Lie algebra ${}^L\g$ is defined so that its root system
is the coroot system $\check\Del$ of $\g$.
In Cartan's classification, ${}^L\g$ is of the same type as $\g$ unless 
$\g$ is of type $B_r$ or $C_r$, in which case ${}^L\g$ is of type $C_r$ or
$B_r$, respectively.
The Lie algebra ${}^L\g$ then determines the universal covering group 
$\widetilde{{}^LG}$, whose centre is 
$Z(\widetilde{{}^LG})\cong\check\Lam^*/\Lam$.
The latter is also the character group 
$Z(\tilde G)^*=\Hom(Z(\tilde G),U(1))$ of $Z(\tilde G)$.
Second, the fundamental group of ${}^LG$, which should be a subgroup of
$Z(\widetilde{{}^LG})\cong Z(\tilde G)^*$, is defined to be
$\pi_1({}^LG)=\set{\chi\in Z(\tilde G)^*\,}{\,\chi(\pi_1(G))=1}$.
$\widetilde{{}^LG}$ and $\pi_1({}^LG)$ uniquely detemine
${}^LG=\widetilde{{}^LG}/\pi_1({}^LG)$.
We have $\pi_1({}^LG)\cong Z(G)^*\cong\ell^*/\Lam$ and
$Z({}^LG)\cong\pi_1(G)^*\cong\check\Lam^*/\ell^*$.
In particular, the Langlands duals of $\tilde G$, $G_\ad$ are $({}^LG)_\ad$,
$\widetilde{{}^LG}$, respectively.

A non-trivial representation of $\g$ is miniscule if all the weights 
form a single Weyl group orbit.
Their highest weights are called miniscule weights.
A miniscule weight must be fundamental, but the converse is not true.
For some simple Lie algebras (of types $B_r$ ($r\ge3$), $D_r$ ($r\ge4$) and
all exceptional types), the adjoint representation is fundamental, but it 
can not be miniscule because $0$ is a weight not in the orbit of the highest
weight.
A fundamental representation may also contain non-zero weights of different
lengths.
For example, the adjoint representation of $B_r$ $(r\ge3$) is fundamental;
it is the second exterior product of the defining representation.
But all non-zero weights of the latter are the weights of the former.
If a representation is miniscule, then so are its conjugate representation
and those obtained by (outer) automorphisms of $\g$.
A complete list of miniscule weights is shown below by highlighting
the corresponding vertices in the Dynkin diagrammes.

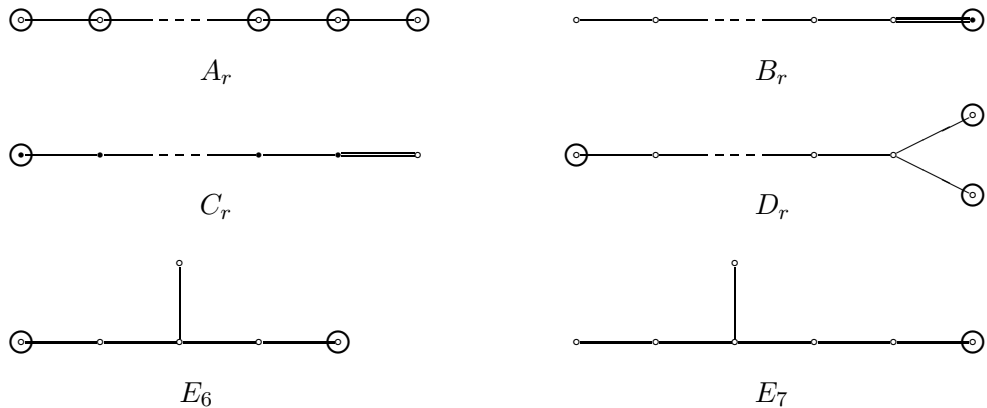
\begin{picture}(100,30)(0,0)
\unitlength 30pt
\put(1.05,0){\line(1,0){.9}}
\multiput(2.05,0)(1.3,0){2}{\line(1,0){.6}}
\multiput(2.75,0)(.2,0){3}{\line(1,0){.1}}
\put(4.05,0){\line(1,0){.9}}
\put(5.05,0){\line(1,0){.9}}
\multiput(1,0)(1,0){2}{\circle{.05}}
\multiput(4,0)(1,0){3}{\circle{.05}}
\put(3.25,-.75){$A_r$}

\put(8.05,0){\line(1,0){.9}}
\multiput(9.05,0)(1.3,0){2}{\line(1,0){.6}}
\multiput(9.75,0)(.2,0){3}{\line(1,0){.1}}
\put(11.05,0){\line(1,0){.9}}
\multiput(12.04,.02)(0,-.04){2}{\line(1,0){.92}}
\put(12.5,0.2){\line(1,-1){.2}}
\put(12.5,-0.2){\line(1,1){.2}}
\multiput(8,0)(1,0){2}{\circle{.05}}
\multiput(11,0)(1,0){2}{\circle{.05}}
\put(13,0){\circle*{.05}}
\put(10.25,-.75){$B_r$}

\thicklines\linethickness{20pt}
\multiput(1,0)(1,0){2}{\circle{.25}}
\multiput(4,0)(1,0){3}{\circle{.25}}
\put(13,0){\circle{.25}}
\end{picture}

\begin{picture}(100,50)(0,0)
\unitlength 30pt
\put(1.05,0){\line(1,0){.9}}
\multiput(2.05,0)(1.3,0){2}{\line(1,0){.6}}
\multiput(2.75,0)(.2,0){3}{\line(1,0){.1}}
\put(4.05,0){\line(1,0){.9}}
\multiput(5.04,.02)(0,-.04){2}{\line(1,0){.92}}
\put(5.5,0.2){\line(-1,-1){.2}}
\put(5.5,-0.2){\line(-1,1){.2}}
\multiput(1,0)(1,0){2}{\circle*{.05}}
\multiput(4,0)(1,0){2}{\circle*{.05}}
\put(6,0){\circle{.05}}
\put(3.25,-.75){$C_r$}

\put(8.05,0){\line(1,0){.9}}
\multiput(9.05,0)(1.3,0){2}{\line(1,0){.6}}
\multiput(9.75,0)(.2,0){3}{\line(1,0){.1}}
\put(11.05,0){\line(1,0){.9}}
\multiput(8,0)(1,0){2}{\circle{.05}}
\multiput(11,0)(1,0){2}{\circle{.05}}
\put(12.04,.02){\line(2,1){.91}}
\put(12.04,-.02){\line(2,-1){.91}}
\multiput(13,.5)(0,-1){2}{\circle{.05}}
\put(10.25,-.75){$D_r$}

\thicklines\linethickness{20pt}
\multiput(1,0)(7,0){2}{\circle{.25}}
\multiput(13,.5)(0,-1){2}{\circle{.25}}
\end{picture}

\begin{picture}(100,70)(0,0)
\unitlength 30pt
\multiput(1.05,0)(1,0){4}{\line(1,0){.9}}
\put(3,.05){\line(0,1){.9}}
\multiput(1,0)(1,0){5}{\circle{.05}}
\put(3,1){\circle{.05}}
\put(3,-.75){$E_6$}

\multiput(8.05,0)(1,0){5}{\line(1,0){.9}}
\put(10,.05){\line(0,1){.9}}
\multiput(8,0)(1,0){6}{\circle{.05}}
\put(10,1){\circle{.05}}
\put(10.25,-.75){$E_7$}

\thicklines\linethickness{20pt}
\multiput(1,0)(4,0){2}{\circle{.25}}
\put(13,0){\circle{.25}}
\end{picture}

\bigskip\bigskip\bigskip

For any simple Lie algebra $\g$, the miniscule weights are in one-to-one 
correspondence with the non-zero elements in $\check\Lam^*/\Lam$
[\ref{B78}, \S VIII.7.3].
To paraphrase this result, notice that the set $\check\Lam^*/\Lam$ is
precisely $Z(\widetilde{{}^LG})$, the centre of (the universal covering of)
the Langlands dual group.
So there is a one-to-one correspondence between the miniscule representations
and the non-trivial elements of $Z(\widetilde{{}^LG})$.
The above list of miniscule representations can be compared with the
table of the centre.

Let $\MM(\g)$ be the set of zero and miniscule weights of $\g$.
The bijection with $Z(\widetilde{{}^LG})=Z(\tilde G)^*$ means that there is
a group structure on the set $\MM(\g)$, which we now describe.
Two irreducible representations of $\g$ are congruent [\ref{D}, \S 3] 
if the differences of their weights are in the root lattice $\Lam$.
This classifies representations of $\g$ according to how $Z(\tilde G)$ acts.
That is, two representations of $\g$ are in the same congruence class if 
and only if they, upon exponentiating to $\tilde G$ and restricting to 
the center, yield the same character. 
(Compare [\ref{D}, Theorem~3.1].)
For example, the representations of $B_r$ fall into two congruence classes:
those that are representations of $SO(2r+1)$ (such as the defining 
representation and its tensor powers) and those that are not (such as 
the spinor representation).
The group structure on the congruence classes comes from the tensor product
(which is well-defined on such classes) and coincides with the multiplication
of characters on $Z(\tilde G)$.
In each congruent class of representations, the (unique) one of the smallest
dimension is the trivial or miniscule representation.
So $\MM(\g)$ is the set of representatives of the congruence classes of
representations.
Moreover, it can shown that given two miniscule representations, their tensor
product contains a miniscule or trivial representation as a summand.

At the group level, we say that a representation of $G$ is miniscule if
the induced representation of the Lie algebra $\g$ is so.
Let $\MM(G)$ be the set of highest weights of the trivial or miniscule
representations of $G$.
Clearly, $\MM(G)$ is a subgroup of $\MM(\g)$ under the tensor product of 
congruence classes, and $\MM(\tilde G)=\MM(\g)$, $\MM(G_\ad)=\{0\}$.
For example, if $G=SU(pq)/\ZA_p$, then $\MM(G)\cong\ZA_q$.
In general, $\MM(G)$ is in one-to-one correspondence with the subset
$\ell^*/\Lam$ of $\check\Lam^*/\Lam\cong Z(\widetilde{{}^LG})$.
Therefore we have the following

\begin{lemma}
There is an isomorphism between the groups $\MM(G)$ and
$\pi_1({}^LG)\cong Z(G)^*$.
\end{lemma}

\sect{Simply laced case: Gauss sum and modular invariance}

As before, let $G$ be a compact, simple Lie group. 
Suppose the Lie algebra $\g$ is simply laced, that is, all the roots are of
the same length.
We consider the Weyl group invariant inner product $(\cdot|\cdot)$ on 
$\ii\gt^*$ so that $(\al|\al)=2$ for all $\al\in\Del$.
We define an analog of the Gauss sum
\be
\GG(G)=\frac{1}{\sqrt{|Z(G)|}}\Sum{\mu\in\MM(G)}\e{\pi\ii(\mu|\mu)}
\ee
and set $\GG(\g)=\GG(\tilde G)$.
For example, $\GG(E_8)=1$ since $Z(E_8)=1$ and
\be
\GG(SU(n))=\frac{1}{\sqrtn}\stackrel{{}_{n-1}}{\Sum{k=0}}
           \e{\pi\ii\frac{k(n-k)}{n}}.
\ee
The latter can be compared with the classical Gauss sum
\be
\stackrel{{}_{n-1}}{\Sum{k=0}}\e{2\pi\ii\frac{k^2}{n}}=\sqrtn\,
                      {\textstyle \frac{1+(-\ii)^n}{1-\ii}}.
\ee

Identifying $\ii\gt$ with $\ii\gt^*$ using the inner product $(\cdot|\cdot)$,
we have $\check\Lam=\Lam\subset\check\Lam^*=\Lam^*$. 
The following is a relation between the Gauss sum of $G$ and that of the
Langlands dual ${}^LG$.

\begin{thm}\label{SIMPLY}
If $\g$ is simply laced and of rank $r$, then
\be\label{dual}
\GG(G)=\overline{\GG({}^LG)}\,\e{\frac{\pi\ii}{4}r}.
\ee
In particular,
\be\label{gauss}
\Sum{\mu\in\MM(\g)}\e{\pi\ii(\mu|\mu)}=\sqrt{|Z(G)|}\,\e{\frac{\pi\ii}{4}r}.
\ee
\end{thm}

\noindent
Examples of (\ref{gauss}) are below.
\begin{center}
\begin{tabular}{|c|c|}\hline
$A_r$ & $\stackrel{{}_r}{\Sum{k=0}}\e{\pi\ii\frac{k(r+1-k)}{r+1}} 
      =\sqrt{r+1}\,\e{\frac{\pi\ii}{4}r}$ \\ \hline
$D_r$ & $1+\e{\pi\ii}+2\e{\pi\ii\frac{r}{4}} 
      =\sqrt{4}\,\e{\frac{\pi\ii}{4}r}$\\ \hline
$E_6$ & $1+2\e{\pi\ii\frac{4}{3}} 
      =\sqrt{3}(\e{\frac{\pi\ii}{4}})^6$ \\ \hline
$E_7$ & $1+\e{\pi\ii\frac{3}{2}}=\sqrt{2}(\e{\frac{\pi\ii}{4}})^7$ \\
      \hline
$E_8$ & $1=(\e{\frac{\pi\ii}{4}})^8$ \\ \hline
\end{tabular}
\end{center}
If $Z(\tilde G)=1$, then (\ref{gauss}) is consistent only if the rank
$r=0\!\!\mod8$, which is the case for $E_8$.
A non-trivial example of (\ref{dual}) is
\be\label{supq}
\GG(p,q)=\sqrt{\frac{q}{p}}\;\e{\frac{\pi\ii}{4}(pq-1)}\,\overline{\GG(q,p)}
\ee
when $G=SU(pq)/\ZA_p$, ${}^LG=SU(pq)/\ZA_q$, $r=pq-1$.
Here
\be
\GG(p,q)=\stackrel{{}_{q-1}}{\Sum{k=1}}\e{\pi\ii p\frac{k(q-k)}{q}}
\ee 
is the generalised Gauss sum [\ref{Ch}, \S V.2].
If $p$, $q$ are distinct odd primes, since 
$\GG(p,q)=\left(\frac{p}{q}\right)\,\GG(1,q)$, where $\left(\frac{p}{q}\right)$
is the Legendre symbol, (\ref{supq}) implies (see for example 
[\ref{Ch}, \S V.3]) the celebrated quadratic reciprocity of Gauss
\be
\left(\frac{p}{q}\right)\,\left(\frac{q}{p}\right)=(-1)^{\frac{(p-1)(q-1)}{4}}.
\ee

The formula (\ref{gauss}) for $G=SU(n)$ appears explicitly in [\ref{VW}];
the generalisation to any simply laced $\g$ is straightforward.
As observed in [\ref{VW}], (\ref{gauss}) is a special case of the following.
Let $V$ be a real vector space with a (not necessarily positive) non-degenerate
quadratic form $(\cdot,\cdot)$, which identifies $V^*$ with $V$.
If $L$ is an even lattice in $V$, then [\ref{vdB}] 
(see also [\ref{MH}, Appen.\ 4])
\be
\Sum{\mu\in L^*/L}\e{\pi\ii(\mu,\mu)}=
\sqrt{|L^*/L|}\;\e{\frac{\pi\ii}{4}\sig},
\ee
where $\sig$ is the signature of $(\cdot,\cdot)$.

(\ref{dual}) can be proved using the Fourier transform on finite Abelian
groups (see for example [\ref{Na}, \S I.4]).
To explain the general setting, let $A$ be a finite Abelian group and
let $A^*=\Hom(A,U(1))$, the group of characters.
The Fourier transform of a function $f\colon A\to\CO$ is 
$\hat f\colon A^*\to\CO$ given by 
$\hat f(\chi)=\sum_{a\in A}\chi(a)f(a)$, $\chi\in A^*$.
The inverse transform is $f=\frac{1}{|A|}\sum_{\chi\in A^*}\hat f(\chi)\chi$.
In particular, $f(0)=\frac{1}{|A|}\sum_{\chi\in A^*}\hat f(\chi)$.
For any subgroup $B$ of $A$, we have the (discrete) Poisson summation formula
[\ref{Na}, Theorem~4.11]
\be\label{dis}
\frac{1}{|B|}\Sum{a\in B}f(a)=\frac{1}{|A|}\Sum{\chi\in(A^*)^B}\hat{f}(\chi),
\ee
where $(A^*)^B\cong(A/B)^*$ is the set of characters of $A$ that are
trivial on $B$.
To apply to our situation, let $A=\Lam^*/\Lam$, $B=\ell/\Lam\cong\MM({}^LG)$.
$A^*$ can be identified with $A$ by $u\mapsto\chi_u$, where
$\chi_u(v)=\e{2\pi\ii(u|v)}$ ($u,v\in A$).
Then $(A^*)^B=\ell^*/\Lam\cong\MM(G)$.
If $f(u)=\e{\pi\ii(u|u)}$ ($u\in A$), then by (\ref{gauss}), for any $v\in A$,
\be
\hat f(\chi_v)=\Sum{u\in\Lam^*/\Lam}\e{\pi\ii(u|u)\tau+2\pi\ii(u|v)}
=|Z(\tilde G)|\,\e{-\pi\ii(v|v)}\,\e{\pi\ii\frac{r}{4}}.
\ee
Hence (\ref{dual}) follows from (\ref{dis}).

We now explain the relation of Theorem~\ref{SIMPLY} with modular invariance.
Recall that the modular group $\Gam=SL(2,\ZA)$ is generated by 
$S={0\;\; -1\choose 1\quad 0}$ and $T={1\quad 1\choose 0\quad 1}$
satisfying the relations
\be\label{rel}
S^2=(ST)^3\in Z(\Gam),\quad S^4=I.
\ee
Thus $\bar\Gam=PSL(2,\ZA)=\ZA_2*\ZA_3$ is the free product of $\ZA_2$
and $\ZA_3$.
We construct a representation of $\Gam$ on the group algebra 
$\CO[\MM(\g)]=\spanc\set{|u\ket\,}{\,u\in\Lam^*/\Lam}$, where $S$ and $T$
act as (using the same notations)
\be\label{action}
\begin{array}{l}
S|u\ket=\frac{1}{\sqrt{|Z(\tilde G)|}}
\Sum{v\in\MM(\g)}\e{-2\pi\ii(u|v)}|v\ket,\\
T|u\ket=\e{-\frac{\pi\ii}{12}r+\pi\ii(u|u)}|u\ket.
\end{array}
\ee
The matrix of $S$ is a discrete Fourier transform, that of $T$ is diagonal,
and both are unitary under the basis $\{|u\ket\}$.
The definition (\ref{action}) is motivated by the S-duality transformations
in the $N=4$ supersymmetric gauge theory [\ref{VW}].
It also coincides with the modular transformations on the characters of
affine Lie algebras at level $1$ [\ref{KP}, \ref{GW}], where $r$ plays the
role of the central charge.

One way to show that (\ref{action}) defines a representation of $\Gam$ 
is to check that the relations in (\ref{rel}) are satisfied.
It is easy to show that $S^2|u\ket=|\!-\!u\ket$ for any $u\in\Lam^*/\Lam$.
Hence $S^2$ commutes with $T$ and $S^4=\id$.
A straightforward calculation yields that, for any $u\in\Lam^*/\Lam$,
\be
\begin{array}{l}
T^{-1}S^{-1}T^{-1}|u\ket=\frac{\e{\frac{\pi\ii}{6}r}}{\sqrt{|Z(\tilde G)|}} 
                         \Sum{v\in\Lam^*/\Lam}\e{-\pi\ii(u-v|u-v)}|v\ket,\\
STS|u\ket=\frac{\e{-\frac{\pi\ii}{12}r}}{|Z(\tilde G)|}\Sum{w\in\Lam^*/\Lam}
          \e{\pi\ii(w|w)}\;\Sum{v\in\Lam^*/\Lam}\e{-\pi\ii(u+v|u+v)}|v\ket.
\end{array}
\ee
So (\ref{action}) satisfies $(ST)^3=S^2$ if and only if (\ref{gauss})
holds [\ref{VW}].

For the Lie group $G$ whose Lie algebra is $\g$, we can define a vector
\be
|G\ket=\frac{1}{\sqrt{|\pi_1(G)|}}\Sum{u\in\pi_1(G)}|u\ket\;
\in\;\CO[\Lam^*/\Lam].
\ee
Then we have the following
\begin{prop}\label{SDUAL}
Under the above notations, we have
\be\label{prop}
S|G\ket=|{}^LG\ket, \quad S|{}^LG\ket=|G\ket.
\ee
\end{prop}
For example, $|\tilde G\ket=|0\ket$ and 
\be
S|0\ket=\frac{1}{\sqrt{|Z(\tilde G)|}}\sum_{u\in\MM(\g)}|u\ket=|G_\ad\ket.
\ee
The proof of (\ref{prop}) is a similar but more general calculation.
\bea
S|G\ket &=& \frac{1}{\sqrt{|\pi_1(G)|}}\sum_{u\in\pi_1(G)}
        \frac{1}{\sqrt{|Z(\tilde G)|}}\sum_{v\in\MM(\g)}
        \e{-2\pi\ii(u|v)}|v\ket \nno
&=& \sqrt{\frac{|\pi_1(G)|}{|Z(\tilde G)|}}\sum_{v\in\pi_1({}^LG)}|v\ket
=|{}^LG\ket.
\eea

To give a more concrete understanding of the consistency of (\ref{action}),
we consider the theta functions [\ref{KP}, \ref{GW}]
\be
\vth_u(z,\tau)=\Sum{\xi\in\Lam+u}\e{\pi\ii(\xi|\xi)\tau+2\pi\ii(\xi|z)},
\ee
where $u\in\Lam^*/\Lam$, $z\in(\ii\gt)^\CO$ and $\tau$ is in the upper 
half plane.
As 
\be
\begin{array}{ll}
\vth_u(z+\lam,\tau)=\e{2\pi\ii(\lam|u)}\,\vth_u(z,\tau),
 & \quad\forall\lam\in\Lam^*, \\
\vth_u(z+\tau\eta,\tau)=\e{-\pi\ii[(\eta|\eta)\tau+2(z|\eta)]}\,\vth_u(z,\tau),
 & \quad\forall\eta\in\Lam,
\end{array}
\ee
the theta functions $\vth_u$ can be considered either as holomorphic
sections of $|Z(\tilde G)|$ line bundles $\LL_u$ (which differ from
each other by flat line bundles) over the Abelian variety 
$A_{\tilde G}=(\ii\gt)^\CO/\Lam^*\oplus\tau\Lam$ or as sections of 
a single line bundle $\LL$, which is the pull-back of each $\LL_u$
from $A_{\tilde G}$ to $A=(\ii\gt)^\CO/\Lam\oplus\tau\Lam$.
Moreover, each $\vth_u$ spans $H^0(A_{\tilde G},\OO(\LL_u))$ and 
$\{\vth_u\}$ is a basis of $H^0(A,\OO(\LL))$.

The modular transformation properties of $\vth_u$ are
\be
\begin{array}{l}
\vth_u(z,\tau+1)=\e{\pi\ii(u|u)}\vth_u(z,\tau),\\
\vth_u\left(\frac{z}{\tau},-\frac{1}{\tau}\right)=
  \frac{1}{\sqrt{|Z(\tilde G)|}}\left(\frac{\tau}{\ii}\right)^{\frac{r}{2}}
  \e{\pi\ii\frac{(z|z)}{\tau}}
  \Sum{v\in\Lam^*/\Lam}\e{-2\pi\ii(u|v)}\,\vth_v(z,\tau);
\end{array}
\ee
the latter is a consequence of the Poisson summation formula.
$\{\vth_u\}$ is an example of vector-valued modular forms [\ref{KM}].
To reproduce the actions of $S$ and $T$ in (\ref{action}), we define
\be
\hat\vth_u(z,\tau,\del)=\e{-2\pi\ii\del}\eta(\tau)^{-r}\vth_u(z,\tau),
\ee
where $\del\in\CO$ and
$\eta(\tau)=\e{\frac{\pi\ii\tau}{12}}\prod_{n=1}^\infty(1-\e{2\pi\ii n\tau})$
is the Dedekind eta-function.
Then
\be
\begin{array}{l}
\hat\vth_u(z,\tau+1,\del)=\e{-\frac{\pi\ii}{12}r+\pi\ii(u|u)}
  \hat\vth_u(z,\tau,\del),\\
\hat\vth_u\left(\frac{z}{\tau},-\frac{1}{\tau},\del+\frac{(z|z)}{2\tau}\right)
  =\frac{1}{\sqrt{|Z(\tilde G)|}}
  \Sum{v\in\Lam^*/\Lam}\e{-2\pi\ii(u|v)}\hat\vth_v(z,\tau,\del),
\end{array}
\ee
in conformity with (\ref{action}).
Notice that the following is an action of $\Gam$ on the triple:
\be
{a\quad b\choose c\quad d}\colon(z,\tau,\del)\mapsto
\left(\frac{z}{c\tau+d},\frac{a\tau+b}{c\tau+d},
\del+\frac{c(z|z)}{2(c\tau+d)}\right).
\ee

We can define, for any compact Lie group $G$ with Lie algebra $\g$,
\bea
\vth_G(z,\tau)&=&\frac{1}{\sqrt{|\pi_1(G)|}}\sum_{u\in\ell/\Lam}\vth_u(z,\tau)
\nno
&=&\frac{1}{\sqrt{|\pi_1(G)|}}
 \sum_{\xi\in\ell}\e{\pi\ii(\xi|\xi)\tau+2\pi\ii(z|\xi)}.
\eea
Then $\vth_G$ is a holomorphic section of a line bundle $\LL_G$ over
$A_G=(\ii\gt)^\CO/(\ell^*\oplus\tau\ell)$.
For example, $\vth_{\tilde G}=\vth_0\in H^0(A_{\tilde G},\OO(\LL_0))$.
In fact, $\vth_G$ spans $H^0(A_G,\OO(\LL_G))$.
By Poisson summation,
\be\label{sdual-th}
\vth_G\left(\frac{z}{\tau},-\frac{1}{\tau}\right)=
\left(\frac{\tau}{\ii}\right)^{\frac{r}{2}}\e{\pi\ii\frac{(z|z)}{\tau}}
\vth_{{}^LG}(z,\tau),
\ee
which reflects (\ref{prop}) in Proposition~\ref{SDUAL}.

With the concrete representation of $\Gam$ on $H^0(A,\OO(\LL))$,
the consistency of (\ref{action}) is automatic and
can be used to prove (\ref{gauss}).
We give a direct proof of (\ref{dual}) from (\ref{sdual-th})
using the technique of M.\ Landsberg (see [\ref{Be}, \S 29]).
Take $\tau=1+\ii\eps$, where $\eps>0$.
Then
\bea
\vth_{{}^LG}(0,\tau)&=&\frac{1}{\sqrt{|\pi_1({}^LG)|}}
   \sum_{\lam\in\ell^*}\e{\pi\ii(\xi|\xi)\tau} \nno
&  =&\frac{1}{\sqrt{|Z(G)|}}\sum_{\mu\in\ell^*/\Lam}\e{\pi\ii(\mu|\mu)}
   \sum_{\xi\in\Lam}\e{-\pi\eps(\xi+\mu|\xi+\mu)}
\eea
because $(\xi|\xi)\in 2\ZA$ and $(\xi|\mu)\in\ZA$.
As $\eps\to0^+$, the sum over $\xi\in\Lam$, multiplied by 
$\eps^{\frac{r}{2}}$, turns to a Gaussian integral over $\xi\in\ii\gt^*$,
whose value is $\frac{1}{\sqrt{\vol(\Lam)}}$.
Therefore to the leading order in $\eps$,
\be
\vth_{{}^LG}(0,\tau)\sim
\frac{1}{\sqrt{|Z(G)|}}\sum_{\mu\in\MM(G)}\e{\pi\ii(\mu|\mu)}
\frac{1}{\eps^{\frac{r}{2}}\sqrt{\vol(\Lam)}}
=\frac{1}{\eps^{\frac{r}{2}}\sqrt{\vol(\Lam)}}\,\GG(G).
\ee
On the other hand, since $-\frac{1}{\tau}\sim-1+\ii\eps$, we have
\be
\vth\left(0,-\frac{1}{\tau}\right)\sim
    \frac{1}{\eps^{\frac{r}{2}}\sqrt{\vol(\Lam)}}\,\overline{\GG({}^LG)},
\ee
and (\ref{dual}) follows.

\sect{Non-simply laced case: the Hecke groups}

For a non-simply laced simple Lie algebra $\g$, let $n_\g$ be the ratio 
of the square lengths of the long and short roots. 
Either $n_\g=2$ (for $B_r$, $C_r$, $F_4$) or $n_\g=3$ (for $G_2$).
Denote by $(\cdot|\cdot)$ the Weyl group invariant inner product on $\ii\gt^*$
such that the square length of the long roots is $2$.
We continue to use it to identify $\ii\gt$ with $\ii\gt^*$.
The root system of the Langlands dual Lie algebra ${}^L\g$ is 
${}^L\!\Del=\frac{1}{\sqrt{n_\g}}\check\Del\subset\ii\gt$.
The (co)root and (co)weight lattices of ${}^L\g$ are
\be
\LLam=\frac{1}{\sqrt{n_\g}}\check\Lam,\quad
(\LLam)\;\check{}=\sqrt{n_\g}\,\Lam,\quad
((\LLam)\,\check{}\,)^*=\frac{1}{\sqrt{n_\g}}\Lam^*,\quad
(\LLam)^*=\sqrt{n_\g}\;\check\Lam^*.
\ee
If the maximal torus of $G$ is $T=\gt/2\pi\ii\ell$ and that of ${}^LG$ is
${}^LT=\gt^*/2\pi\ii{}^L\ell$, then ${}^L\ell=\sqrt{n_\g}\,\ell^*$.
We have the inclusion relations 
\be
(\LLam)\,\check{}\subset{}^L\ell\subset(\LLam)^*\subset\ii\gt^*,\quad
\LLam\subset({}^L\ell)^*\subset((\LLam)\,\check{}\,)^*\subset\ii\gt.
\ee

We make a digression on Coxeter and dual Coxeter numbers.
For every simple Lie algebra $\g$, the Coxeter number $h$ and the dual
Coxeter number $\check h$ of $\g$ are listed as follows.
\begin{center}
\begin{tabular}{|c|ccccccccc|}\hline
$\g$ & $A_r$ & $B_r$ & $C_r$ & $D_r$ & $E_6$ & $E_7$ & $E_8$ & $F_4$ & $G_2$ \\
\hline
$h$ & $r+1$ & $2r$ & $2r$ & $r+1$ & $12$ & $18$ & $30$ & $12$ & $6$ \\ \hline  
$\check h$ & $r+1$ & $2r-1$ & $r+1$ & $2r-2$ & $12$ & $18$ & $30$ & $9$ & $4$
\\ \hline
\end{tabular}
\end{center}
We have $\check h=1+\bra\rho,\check\theta\ket$, where $\rho$ is half
the sum of positive roots and $\theta$, the highest root.
If $\g$ is simply laced, then $h=\check h$.
If $\g$ is non-simply laced, let $\rho\loong$, $\rho\short$ be half the sums
of long, short positive roots, respectively, and let 
$h\loong=1+\bra\rho\loong,\check\theta\ket$,
$h\short=n_\g\,\bra\rho\short,\check\theta\ket$.
Then $h=h\loong+h\short$ and $\check h=h\loong+n_\g^{-1}h\short$.
When $\Del$ is exchanged with $\check\Del$, or $\g$ with ${}^L\g$,
$h\loong$ and $h\short$ are also exchanged [\ref{S}].
Thus the Coxeter number ${}^Lh$ and the dual Coxeter number ${}^L\check h$
of the Langlands dual ${}^L\g$ satisfy 
\be\label{cox}
{}^Lh=h, \quad \check h+{}^L\check h=(1+n_\g^{-1})\,h.
\ee
It is well known [\ref{B456}, Th\'eor\`eme~V.6.2.1(ii), Exer.\ VI.1.20] that
$|\Del|=r\,h$, $|\Del\loong|=r\loong\,h$, $|\Del\short|=r\short\,h$, where
$\Del\loong$, $\Del\short$ are the sets of long, short roots and $r\loong$,
$r\short$ are the numbers of long, short simple roots, respectively.
In addition, we have $|\Del\loong|=r\,h\loong$, $|\Del\short|=r\,h\short$
[\ref{S}].

To find the counterpart of the modular group and the correct analog of
(\ref{action}) in the non-simply laced case, we proceed with the
theta functions.
In fact, we need two sets of such.
For each $u\in\Lam^*/\check\Lam$ and $\mu\in(\LLam)^*/(\LLam)\,\check{}\,$, let
\be
\begin{array}{l}
\vth_u(z,\tau)=\Sum{x\in\Lam+u}\e{\pi\ii(x|x)\tau+2\pi\ii(x|z)},\\
\vth_\mu(z,\tau)=\Sum{\lam\in(\LLam)\,\check{}\,+\mu}
   \e{\pi\ii(\lam|\lam)\tau+2\pi\ii(\lam|z)}.
\end{array}
\ee
Here $z\in(\ii\gt)^\CO\cong(\ii\gt^*)^\CO$, where the identification is
via $(\cdot|\cdot)$.
As in the simply laced case, the periodic properties
\be
\begin{array}{ll}
\vth_u(z+\lam,\tau)=\e{2\pi\ii(\lam|u)}\vth(z,\tau),
& \quad\forall\lam\in\check\Lam^*, \\
\vth_u(z+\tau y,\tau)=\e{-\pi\ii[(y|y)\tau+2(z|y)]}\,\vth_u(z,\tau),
& \quad\forall y\in\check\Lam
\end{array}
\ee
mean that for each $u\in\Lam^*/\check\Lam$, $\vth_u$ spans 
$H^0(A_0,\OO(\LL_u))$, where $\LL_u$ is a line bundle over 
$A_{\tilde G}=(\ii\gt)^\CO/\check\Lam^*\oplus\check\Lam\tau$ and that
$\{\vth_u\}$ is a basis of $H^0(A,\OO(\LL))$, where $\LL$ is the pull-back
of (any) $\LL_u$ to $A=(\ii\gt)^\CO/\Lam\oplus\check\Lam\tau$.
Similarly, for each $\mu\in(\LLam)^*/(\LLam)\,\check{}\,$, $\vth_\mu$
spans $H^0(\LA_{\tilde G},\OO(\LLL_\mu))$ and $\{\vth_\mu\}$ is 
a basis of $H^0(\LA,\OO(\LLL))$.
Here $\LLL\to\LA=(\ii\gt)^\CO/\LLam\oplus\tau(\LLam)\,\check{}\,$
is the pull-back of the line bundle $\LLL_\mu\to\LA_{\tilde G}=
(\ii\gt)^\CO/((\LLam)\,\check{}\,)^*\oplus\tau(\LLam)\,\check{}\,$.

It is easy to see that, using $(\check\Lam|\check\Lam)\subset 2\ZA$,
$(\Lam^*|\check\Lam)\subset\ZA$ and similar facts for $\LLam$, for any
$u\in\Lam^*/\check\Lam$ and $\mu\in(\LLam)^*/(\LLam)\,\check{}\,$,
\be
\begin{array}{l}
\vth_u(z,\tau+1)=\e{\pi\ii(u|u)}\vth_u(z,\tau),\\
\vth_\mu(z,\tau+1)=\e{\pi\ii(\mu|\mu)}\vth_\mu(z,\tau).
\end{array}
\ee
Poisson summation yields
\be
\begin{array}{l}
\vth_u\left(\frac{z}{\sqrt{n_\g}\tau},-\frac{1}{n_\g\tau}\right) 
  =\sqrt{\frac{n_\g^{r\loong}}{|Z(\tilde G)|}}
  \left(\frac{\tau}{\ii}\right)^{\frac{r}{2}}
  \e{\pi\ii\frac{(z|z)}{\tau}}\Sum{\mu\in(\LLam)^*/(\LLam)\;\check{}}
  \e{-2\pi\ii\frac{\bra\mu,u\ket}{\sqrt{n_\g}}}\vth_\mu(z,\tau), \\
\vth_\mu\left(\frac{z}{\sqrt{n_\g}\tau},-\frac{1}{n_\g\tau}\right)
  =\sqrt{\frac{n_\g^{r\short}}{|Z(\widetilde{{}^LG})|}}
  \left(\frac{\tau}{\ii}\right)^{\frac{r}{2}}
  \e{\pi\ii\frac{(z|z)}{\tau}}\Sum{u\in\Lam^*/\check\Lam}
  \e{-2\pi\ii\frac{\bra\mu,u\ket}{\sqrt{n_\g}}}\vth_u(z,\tau).
\end{array}
\ee
Here the factor $n_\g$ appears in order to convert the sum over $\check\Lam^*$
to $(\LLam)^*=\sqrt{n_\g}\,\check\Lam^*$.
We define
\be
\begin{array}{l}
\hat\vth_u(z,\tau,\del)=\e{-2\pi\ii\del}\,\eta(\tau)^{-r\loong}\,
  \eta(n_g\tau)^{-r\short}\,\vth_u(z,\tau),\\
\hat\vth_\mu(z,\tau,\del)=\e{-2\pi\ii\del}\,\eta(\tau)^{-r\short}\,
  \eta(n_g\tau)^{-r\loong}\,\vth_\mu(z,\tau).
\end{array}
\ee
Then
\be\label{he-th}
\begin{array}{l}
\hat\vth_u(z,\tau+1,\del)=
  \e{-\frac{\pi\ii}{12}n_\g\frac{{}^L\check h}{h}r+\pi\ii(u|u)}
  \hat\vth_u(z,\tau,\del),\\
\hat\vth_\mu(z,\tau+1,\del)=
  \e{-\frac{\pi\ii}{12}n_\g\frac{\check h}{h}r+\pi\ii(\mu|\mu)}
  \hat\vth_\mu(z,\tau,\del),\\
\hat\vth_u\left(\frac{z}{\sqrt{n_\g}\tau},-\frac{1}{n_\g\tau},
  \del+\frac{(z|z)}{2\tau}\right)
  =\frac{1}{\sqrt{|Z(\tilde G)|}}\Sum{\mu\in(\LLam)^*/(\LLam)\;\check{}}
  \e{-2\pi\ii\frac{\bra\mu,u\ket}{\sqrt{n_\g}}}\hat\vth_\mu(z,\tau,\del),\\
\hat\vth_\mu\left(\frac{z}{\sqrt{n_\g}\tau},-\frac{1}{n_\g\tau},
  \del+\frac{(z|z)}{2\tau}\right)
  =\frac{1}{\sqrt{|Z(\tilde G)|}}\Sum{u\in\Lam^*/\check\Lam}
  \e{-2\pi\ii\frac{\bra\mu,u\ket}{\sqrt{n_\g}}}\hat\vth_u(z,\tau,\del).
\end{array}
\ee

We can also define the theta functions
\be
\begin{array}{rl}
\vth_G(z,\tau)&=
 \frac{1}{\sqrt{|\pi_1(G)|}}\Sum{u\in\ell/\check\Lam}\vth_u(z,\tau)\\
&=\frac{1}{\sqrt{|\pi_1(G)|}}\Sum{x\in\ell}\e{\pi\ii(x|x)\tau+2\pi\ii(z|x)},\\
\vth_{{}^LG}(z,\tau)&=\frac{1}{\sqrt{|\pi_1({}^LG)|}}
   \Sum{\mu\in{}^L\ell/(\LLam)\;\check{}}\vth_\mu(z,\tau) \\
&=\frac{1}{\sqrt{|\pi_1({}^LG)|}}
   \Sum{\lam\in{}^L\ell}\e{\pi\ii(\lam|\lam)\tau+2\pi\ii(z|\lam)}.\\
\end{array}
\ee
Poisson summation yields
\be
\vth_G\left(\frac{z}{\sqrt{n_\g}\tau},-\frac{1}{n_\g\tau}\right)=
\sqrt{n_\g^{r\loong}}\left(\frac{\tau}{\ii}\right)^{\frac{r}{2}}
\e{\pi\ii\frac{(z|z)}{\tau}}\vth_{{}^LG}(z,\tau).
\ee
$\vth_G$ is a holomorphic section of a line bundle 
$\LL_G\to A_G=(\ii\gt)^\CO/\ell^*\oplus\tau\ell$
while $\vth_{{}^LG}$ is that of a line bundle 
$\LL_{{}^LG}\to A_{{}^LG}=(\ii\gt)^\CO/({}^L\ell)^*\oplus\tau{}^L\ell$.
It is curious to observe that, for both simply laced and non-simply laced
cases,  the Abelian variety $A_{{}^LG}$ for $\tau$ is dual to $A_G$ for
$-\frac{1}{n_\g\tau}$ and that $\LL_{{}^LG}$ is the Fourier-Mukai
transform [\ref{M}] of $\LL_G$.  

The transformations $\tau\mapsto\tau+1$, $\tau\mapsto-\frac{1}{n_\g\tau}$
generate what is called the Hecke group [\ref{He}].
Recall that the Hecke group $G(\lam_m)$ ($3\le m\le+\infty$) is the subgroup 
of $SL(2,\RE)$ generated by $S$ and 
$\tilde T={1\,\;\;\lam_m\choose 0\quad 1\;}$,
where $\lam_n=2\cos\frac{\pi}{m}$, satisfying the relations
\be
S^2=(S\tilde T)^m\in Z(G(\lam_m)),\quad S^4=1.
\ee
The image $\bar G(\lam_m)$ in $PSL(2,\RE)$ is the free product of $\ZA_2$
and $\ZA_m$.
When $m=3$, we have the classical modular group $\Gam=SL(2,\ZA)$.
When $m=4,6$, the groups are $G(\sqrt2)$, $G(\sqrt3)$, respectively.
Consider also the subgroup $\Gam^*(n)$ of $SL(2,\RE)$ generated by
$\tilde S={\;\;0\;\;-1/\sqrtn\choose\!\!\sqrtn\;\quad\;0}$ and $T$, 
which define the desired fractional linear transformations 
$\tau\mapsto\tau+1$, $\tau\mapsto-\frac{1}{n\tau}$ on $\tau$.
The group $\Gam^*(n)$ contains 
$\Gam_0(n)=\set{{a\quad b\choose c\quad d}}{\,c=0\!\!\mod n}$ as a subgroup 
of index $2$.
When $n$ is prime, $\Gam^*(n)$ is the normaliser of $\Gam_0(n)$ in $SL(2,\RE)$ 
[\ref{LN}].
When $n=n_\g=2,3$, $\Gam^*(n_g)$ is isomorphic to $G(\sqrt{n_g})$ by a 
conjugation with the matrix ${1\quad\;0\;\choose 0\;\sqrt{n_g}}$, under 
which the generators $S$, $\tilde T$ of $G(\sqrt{n_\g})$ are mapped to 
$\tilde S$ and $T$, respectively.
Hence $\tilde S$ and $T$ satisfy the same relations
\be\label{hecke}
\tilde S^2=(\tilde ST)^{2n_\g}\in Z(\Gam^*(n_g)),\quad \tilde S^4=1.
\ee
These two groups, together with $G(1)=\Gam$ and $G(2)=\Gam(2)$, are the only
Hecke groups $G(\lam_m)$ that are commensurable with the modular group 
[\ref{Le}].

By the transformations (\ref{he-th}), the Hecke group $\Gam^*(n_\g)$ acts 
on the direct sum $H^0(A,\OO(\LL))\oplus H^0(\LA,\OO(\LLL))$.
The subgroup $\Gam_0(n_\g)$ of index $2$ acts within the spaces 
$H^0(A,\OO(\LL))$ and $H^0(\LA,\OO(\LLL))$ while $\tilde S$ exchanges them.
More abstractly, we define a representation of $\Gam^*(n_\g)$ on\\
$\CO[\Lam^*/\check\Lam]\oplus\CO[(\LLam)^*/(\LLam)\,\check{}\;]
=\spanc\set{|u\ket,|\mu\ket\,}{\,u\in\Lam^*/\check\Lam,
\mu\in(\LLam)^*/(\LLam)\,\check{}\,}$
such that $\tilde S$ and $T$ act as (using the same notations)
\be\label{he-repr}
\begin{array}{l}
T|u\ket=\e{-\frac{\pi\ii}{12}n_\g\frac{{}^L\check h}{h}r+\pi\ii(u|u)}|u\ket,\\
T|\mu\ket=
\e{-\frac{\pi\ii}{12}n_\g\frac{\check h}{h}r+\pi\ii(\mu|\mu)}|\mu\ket,\\
\tilde S|u\ket=\frac{1}{\sqrt{|Z(\tilde G)|}}
  \Sum{\mu\in(\LLam)^*/(\LLam)\,\check{}}
  \e{-2\pi\ii\frac{\bra\mu,u\ket}{\sqrt{n_\g}}}|\mu\ket,\\
\tilde S|\mu\ket=\frac{1}{\sqrt{|Z(\tilde G)|}}\Sum{u\in\Lam^*/\check\Lam}
  \e{-2\pi\ii\frac{\bra\mu,u\ket}{\sqrt{n_\g}}}|u\ket.
\end{array}
\ee
It is easy to see that $\tilde S$ and $T$ are unitary and that
$\tilde S^2|u\ket=|\!-\!u\ket$, $\tilde S^2|\mu\ket=|\!-\!\mu\ket$.
Because of the concrete realisation by the theta functions, the above formulae
are consistent with the relations (\ref{hecke}).
We can check $(\tilde ST)^{2n_\g}=\id$ explicitly in all non-simply laced
cases.
Since the centre is either $\ZA_2$ or trivial, the action of $S^2$ is
the identity.
The contribution of the new phase factors 
$\e{-\frac{\pi\ii}{12}n_\g\frac{{}^L\check h}{h}r}$ and
$\e{-\frac{\pi\ii}{12}n_\g\frac{\check h}{h}r}$ to $(\tilde ST)^{2n_\g}$ is
\be
\left(\e{-\frac{\pi\ii}{12}n_\g\frac{{}^L\check h}{h}r}\,
\e{-\frac{\pi\ii}{12}n_\g\frac{\check h}{h}r}\right)^{n_\g}=
\e{-\frac{\pi\ii}{12}n_\g(n_\g+1)r}
\ee
by using (\ref{cox}).
Thus $(\tilde ST)^{2n_\g}$ can be calculated as follows.
\begin{center}
\begin{tabular}{|c|c|}\hline
$B_r$ or $C_r$ & 
$(\tilde ST)^4=\e{-\frac{\pi\ii}{12}2\cdot3\cdot r}{\e{\pi\ii\frac{r}{2}}
\quad\quad\quad\choose\quad\quad\quad\e{\pi\ii\frac{r}{2}}}
={1\;\quad{}\choose{}\;\quad1}$ \\ \hline
$F_4$ & $(\tilde ST)^4=\e{-\frac{\pi\ii}{12}2\cdot3\cdot4}=1$ \\ \hline
$G_2$ & $(\tilde ST)^6=\e{-\frac{\pi\ii}{12}3\cdot4\cdot2}=1$ \\ \hline
\end{tabular}
\end{center}
If we define
\be\label{glg}
\begin{array}{l}
|G\ket=\frac{1}{\sqrt{|\pi_1(G)|}}\Sum{u\in\pi_1(G)}|u\ket
  \in\CO[\Lam^*/\check\Lam], \\
|{}^LG\ket=\frac{1}{\sqrt{|\pi_1({}^LG)|}}\Sum{\mu\in\pi_1({}^LG)}
  |\mu\ket\in\CO[(\LLam)^*/(\LLam)\,\check{}\,],
\end{array}
\ee
then $\tilde S|G\ket=|{}^LG\ket$ and $\tilde S|{}^LG\ket=G\ket$ as before.

We summarize the results in the following
\begin{thm}
If $\g$ is non-simply laced, (\ref{he-repr}) defines a unitary
representation of the Hecke group $\Gam^*(n_\g)\cong G(\sqrt{n_\g})$ on 
$\CO[\Lam^*/\check\Lam]\oplus\CO[(\LLam)^*/(\LLam)\,\check{}\,]$.
Moreover $\tilde S$ interchanges $|G\ket$ and $|{}^LG\ket$ in (\ref{glg}).
\end{thm}

Finally, we note some new features of the non-simply laced case.
The transformations (\ref{he-th}) of the theta functions do not coincide
with those of the character of affine Lie algebras.
It would be interesting to find an alternative relation and to establish
an analog of the vector-valued modular forms [\ref{KM}] for Hecke groups.
Unlike the simply laced case, the new phase factors in the transformation
of $T$ do not come from the central charges of conformal field theory.
In a forthcoming work, we will explore the consequences of this in $S$-duality.

\bigskip

\noindent
{\small {\bf Acknowledgements.}
The author would like to thank X.\ Wang and Y.\ Zhu for discussions.
This work is supported in part by a CERG grant HKU705407P.}

        \newcommand{\athr}[2]{{#1}.\ {#2}}
        \newcommand{\au}[2]{\athr{{#1}}{{#2}},}
        \newcommand{\an}[2]{\athr{{#1}}{{#2}} and}
        \newcommand{\jr}[6]{{#1}, {\it {#2}} {#3}\ ({#4}) {#5}-{#6}}
        \newcommand{\pr}[3]{{#1}, {#2} ({#3})}
        \newcommand{\bk}[4]{{\it {#1}}, {#2}, {#3} ({#4})}
        \newcommand{\cf}[8]{{\it {#1}}, {#2}, {#5},
                 {#6}, ({#7}, {#8}), pp.\ {#3}-{#4}}

        \vspace{3ex}
        \begin{flushleft}
{\bf References}
        \end{flushleft}
{\small
        \baselineskip=11pt
\begin{enumerate}

\item\label{Be}
\au{R}{Bellman}
\bk{A brief introduction to theta functions}
{Holt, Rinehart and Winston}{New York}{1961}.

\item\label{vdB}
\au{F}{van der Blij}
\jr{An invariant of quadratic forms mod $8$}
{Indag.\ Math.}{21}{1959}{291}{293}.

\item\label{B456}
\au{N}{Bourbaki}
\bk{Groupes et alg\`ebres de Lie, Chap.\ IV, V et VI}
{Hermann}{Paris}{1968}.

\item\label{B78}
\au{N}{Bourbaki}
\bk{Groupes et alg\`ebres de Lie, Chap.\ VII et VIII}
{Hermann}{Paris}{1975}.

\item\label{B9}
\au{N}{Bourbaki}
\bk{Groupes et alg\`ebres de Lie, Chap.\ IX}{Masson}{Paris}{1982}

\item\label{Ch}
\au{K}{Chandrasekharan}
\bk{Introduction to analytic number theory}
{Springer-Verlag}{Berlin-New York}{1968}.

\item\label{D}
\au{E}{Dynkin}
\jr{Semisimple subalgebras of semisimple Lie algebras}
{Mat.\ Sbornik N.S.}{30(72)}{1952}{349}{462};
\jr{English translation}{Amer.\ Math.\ Soc.\ Transl.\ (Ser.\ 2)}
{6}{1957}{111}{244}.

\item\label{GW}
\an{D}{Gepner} \au{E}{Witten}
\jr{String theory on group manifolds}{Nucl.\ Phys.\ B}{278}{1986}{493}{549}.

\item\label{GNO}
\au{P}{Goddard} \an{J}{Nuyts} \au{D.\ I}{Olive}
\jr{Gauge theories and magnetic charges}{Nucl.\ Phys.\ B}{125}{1977}{1}{28}.

\item\label{He}
\au{E}{Hecke}
\jr{\"Uber die Bestimmung Dirichletscher Reihen durch ihre Funktionalgleichung}
{Math.\ Ann.}{112}{1936}{664}{699}.

\item\label{KP}
\an{V.\ G}{Ka\v c} \au{D.\ H}{Peterson}
\jr{Infinite-dimensional Lie algebras, theta functions and modular forms}
{Adv.\ Math.}{53}{1984}{125}{264}.

\item\label{KM}
\an{M}{Knopp} \au{G}{Mason}
\jr{Vector-valued modular forms and Poincar\'e series}
{Illinois J.\ Math.}{48}{2004}{1345}{1366}.

\item\label{LN}
\an{J}{Lehner} \au{M}{Newman}
\jr{Weierstrass points of $\Gam_0(n)$}{Ann.\ Math.\ (2)}{79}{1964}{360}{368}.

\item\label{Le}
\au{A}{Leutbecher}
\jr{\"Uber die Heckeschen Gruppen $G(\lambda)$}
{Abh.\ Math.\ Sem.\ Univ.\ Hamburg}{31}{1967}{199}{205};
\jr{Teil II}{Math.\ Ann.}{211}{1974}{63}{86}.

\item\label{MH}
\an{J}{Milnor} \au{D}{Husemoller}
\bk{Symmetric bilinear forms}{Springer-Verlag}{New York-Heidelberg}{1973}.

\item\label{MO}
\an{C}{Montonen} \au{D.\ I}{Olive}
\jr{Magnetic monopoles as gauge particles?}
{Phys.\ Lett.\ B}{72}{1977}{117}{120}.

\item\label{M}
\au{S}{Mukai}
\jr{Duality between $D(X)$ and $D(\hat X)$ with its application to Picard
sheaves}{Nagoya Math.\ J.}{81}{1981}{153}{175}.

\item\label{Na}
\au{M.\ B}{Nathanson}
\bk{Elementary methods in number theory}{Springer}{New York}{2000}.

\item\label{S}
\au{R}{Suter}
\jr{Coxeter and dual Coxeter numbers}{Commun.\ Alg.}{26}{1998}{147}{153}.

\item\label{VW}
\an{C}{Vafa} \au{E}{Witten}
\jr{A strong coupling test of $S$-duality} 
{Nucl.\ Phys.\ B}{431}{1994}{3}{77}.

\end{enumerate}}

\end{document}